\documentclass{article}
\usepackage{amsfonts}
\usepackage{amsmath}
\usepackage{amssymb}
\usepackage{epsfig}
\usepackage{verbatim}   
\usepackage{color}      
\usepackage{hyperref}   

\begin{document}

\title{Noncommutative Biorthogonal Polynomials}
\author{Emily Sergel
\footnote{Special thanks to Prof. Robert Wilson, Rutgers University, for overseeing and guiding this research, as well as all assistance in creating this paper. This research was funded by NSF-0603745.} \\Rutgers University}

\maketitle

\vskip 12 pt
The theory of orthogonal polynomials is well established and has many applications. For any sequence $\{S_i\}$ of elements of a commutative ring $R$, we can define a biadditive function $< \cdot, \cdot > : R[x] \times R[x] \rightarrow R$ by $<ax^i,bx^j>=abS_{i+j}$ for $a,b \in R$ and define a sequence of polynomials $\{p_n\}$ by $$p_{n} = \vmatrix S_{n} & \cdots & S_{2n-1} & x^{n} \\ \vdots & \ddots & \vdots & \vdots \\ S_{0} & \cdots & S_{n-1} & 1 \endvmatrix.$$ Then $<p_n,p_m.$ if and only if $n \neq m$, i.e. the sequence $\{p_n\}$ is orthogonal. The $S_i$ are called the moments of $\{p_n\}$. For a more detailed introduction, see Chihara's classic text {\cite{C}} on the subject. The idea of orthogonal polynomials and this method of generating them has been generalized in two ways to achieve new types of polynomials: noncommutative orthogonal polynomials and biorthogonal polynomials.\\

The theory of orthogonal polynomials has been extended to cover rings of noncommutative operators, such as matrices (see {\cite{M}}). In 1994, Gelfand, Krob, Lascoux, Leclerc, Retakh and Thibon{\cite{GKLLRT}} extended the theory to general noncommutative rings by setting $p_n$ equal to the quasideterminant of a similar matrix. The paper also shows that the 3-term recurrence relation, which is well-known for commutative orthogonal polynomials, still holds in this case.\\

Second, orthogonal polynomials have been generalized in several ways to biorthogonal polynomials.  See {\cite{B}} for more details on these generalizations. One such extension is considered by Bertola, Gekhtman and Szmigielski{\cite{BGS}}. A family of biorthogonal polynomials is defined to be two sequences of real polynomials $\{p_{n}(x)\}$ and $\{q_{m}(y)\}$ with the property that $\int\int p_{n}(x)q_{m}(y)K(x,y) d\alpha(x)d\beta(y) = 0$ when $n \neq m$ for particular K, $\alpha$ and $\beta$.  In this paper, it is shown that these polynomials can be represented as determinants of matrices whose entries are bimoments and, for a specific $K(x,y)$, a 4-term recurrence relation is obtained.\\

Here, we define biorthogonal polynomials over a noncommutative ring. We bring together the two different generalizations described above to present a completely algebraic definition of noncommutative biorthogonal polynomials. For our purposes, a biorthogonal family consists of two sequences of polynomials $\{p_n(x)\}$ and $\{q_n(y)\}$, over a division ring R, along with a function $<\cdot,\cdot>$: $R[x] \times R[y] \rightarrow R$ so that $<p_n(x),q_m(y)> = 0$ for all $n \neq m$. Using this definition, we obtain recurrence relations for some types of biorthogonal polynomials and thus generalize the 4-term recurrence relations of {\cite{BGS}}. We conclude with a broad extension of Favard's theorem.
\pagebreak

\section{Set-Up and Definitions}

Let R be a division ring with center C. We will view R[x] as an R-C bimodule of R and R[y] as a C-R bimodule of R. That is, elements of R[x] will be of the form $\sum a_{i}x^{i}$ and elements of R[y] will be of the form $\sum y^{j}b_{j}$ so that $xc=cx$ and $yc=cy$ for all $c \in C$. Let $<\cdot,\cdot> : R[x] \times R[y] \rightarrow R$ so that $$<\sum a_{i}x^{i}, \sum y^{j}b_{j}> = \sum a_{i}<x^{i},y^{j}>b_{j}.$$
A system of polynomials $\{p_{n}\}, \{q_{n}\}_{n \in \mathbb{N}}$ is \underline{biorthogonal} with respect to $<\cdot,\cdot>$ if $<p_{n}(x), q_{m}(y)> = 0$ for all $n \neq m.$\\
\vskip 12 pt
Let $I_{a,b} = <x^{a}, y^{b}>.$ The set $I = \{I_{a,b}\}_{a,b \in \mathbb{Z}_{\geq 0}}$ is called the set of \underline{bimoments} for $<\cdot,\cdot>.$ The bimoments completely define the function $<\cdot, \cdot>$ so we will say that a set of polynomials is biorthogonal with respect to $I$.\\
In keeping with the notation of {\cite{BGS}}, we will let $I$ be the matrix of bimoments and write $Id$ for the identity matrix. Note in these cases, and below, all matrices and vectors are infinite, with rows and columns indexed by $\mathbb{Z}_{\geq 0}$.\\\\
\vskip 12 pt
We extend $<\cdot,\cdot>$ to $R[x]^{n} \times R[y]$ and to $R[x] \times R[y]^n$ in the following way:\\
If $B = \bmatrix b_{1} \\ \vdots \\ b_{n} \endbmatrix \in R[x]^n $ and $g \in R[y], then <B,g> = \bmatrix <b_{1},g> \\ \vdots \\ <b_{n},g> \endbmatrix$.\\
\vskip 6 pt
Similarly, if $f \in R[x]$ and $D = \bmatrix d_{1} \\ \vdots \\ d_{n} \endbmatrix \in R[y]^n$, then $<f,D> = \bmatrix <f,d_{1}> \\ \vdots \\ <f,d_{n}> \endbmatrix$.\\
\vskip 6 pt
If $C \in Mat_{r \times n}(R), B \in R[x]^n$ and $g \in R[y]$, then $<CB,g> = C<B,g>$.\\
\vskip 12 pt
For an (n+1)x(n+1) matrix A, let $A^{i,j}$ denote the nxn matrix formed by removing the ith row and jth column. Then(c.f. {\cite{GKLLRT}}, def. 2.1) the \underline{i,j-quasideterminant of A} $|A|_{i,j}$ is $$\vmatrix a_{1,1} & \cdots & a_{1,j} & \cdots & a_{1,n+1}\\ \vdots & \ddots & \vdots & \ddots & \vdots \\ a_{i,1} & \cdots & $\fbox{$a_{i,j}$}$ & \cdots & a_{i,n+1}\\ \vdots & \ddots & \vdots & \ddots & \vdots \\ a_{n+1,1} & \cdots & a_{n+1,j} & \cdots & a_{n+1,n+1} \endvmatrix$$ $$= a_{i,j} - \bmatrix a_{i,1} & \cdots & a_{i,j-1} & a_{i,j+1} & \cdots & a_{i,n+1} \endbmatrix \cdot (A^{i,j})^{-1} \cdot \bmatrix a_{1,j} \\ \cdots \\ a_{i-1,j} \\ a_{i+1,j} \\ \cdots \\ a_{n+1,j} \endbmatrix.$$
Note that, after suitably permuting rows and columns, this the is Schur complement of a block decomposition of $A$. The quasideterminant $|A|_{i,j}$ exists if and only if $A^{i,j}$ is invertible.
\vskip 12 pt

\section{Constructing Biorthogonal Polynomials Using Bimoments}
Throughout, we will assume that the set of bimoments is generic in the sense that all quasideterminants considered exist and are invertible. This is our only restriction on the set of bimoments.

\textbf{Theorem:}
Let $\{I_{a,b}|a,b \in \mathbb{Z}_{\geq 0}\} \subseteq R.$
For all $n \in \mathbb{N}$, define $$p_{n}(x) = |I|_{1,n+1} = \vmatrix I_{n,0} & \cdots & I_{n,n-1} & $\fbox{$x^{n}$} $\\ \vdots & \ddots & \vdots & \vdots \\ I_{1,0} & \cdots & I_{1,n-1} & x \\ I_{0,0} & \cdots & I_{0,n-1} & 1 \endvmatrix
$$ and $$ q_{n}(y) = \vmatrix 1 & y & \cdots & $\fbox{$y^{n}$}$ \\ I_{n-1,0} & I_{n-1,1} & \cdots & I_{n-1,n} \\ \vdots & \vdots & \ddots & \vdots \\ I_{0,0} & I_{0,1} & \cdots & I_{0,n} \endvmatrix.$$ Then $\{p_{n}\}, \{q_{n}\}$ is a (monic) biorthogonal system of polynomials with respect to the set of bimoments $\{I_{a,b}\}$.\\

To prove the theorem we need the following lemma:
\vskip 12 pt
\textbf{Lemma:}
Let $n \in \mathbb{Z}_{\geq 0}$ and $p_{n}, q_{n}$ be as defined as in the proposition. Then $<x^{i},q_{n}> = <p_{n},y^{i}> = 0$ for all $0 \leq i \leq n-1$.

\vskip 12 pt
\textbf{Proof of Lemma:}
\\
Let $n \in \mathbb{N}$ and $0 \leq i \leq n-1$.
We see that $$<p_{n},y^{i}>  =  <x^{n} - \bmatrix I_{n,0} & \cdots & I_{n,n-1} \endbmatrix \cdot (I^{1,n+1})^{-1} \cdot \bmatrix x^{n-1} \\ \vdots \\ 1 \endbmatrix, y^{i}> $$ = $$I_{n,i} - \bmatrix I_{n,0} & \cdots & I_{n,n-1} \endbmatrix \cdot (I^{1,n+1})^{-1} \cdot \bmatrix I_{n-1,i} \\ \vdots \\ I_{0,i} \endbmatrix.$$ Applying the definition of quasideterminant, we see that this is $$ \vmatrix I_{n,0} & \cdots & I_{n,n-1} & $\fbox{$I_{n,i}$}$ \\ \vdots & \ddots & \vdots & \vdots \\ I_{0,0} & \cdots & I_{0,n-1} & I_{0,i} \endvmatrix.$$
Thus, since $0 \leq i \leq n-1$, $<p_{n}, y^{i}>$ is the quasideterminant of a matrix whose nth column is equal to its $(i+1)st$ column and hence is 0 (c.f. {\cite{GGRW}}, prop. 1.4.6).\\
Similarly, $$<x^{i}, q_{n}> = \vmatrix I_{i,0} & \cdots & $\fbox{$I_{i,n}$} $\\ I_{n-1,0} & \cdots & I_{n-1,n} \\ \vdots & \ddots & \vdots \\ I_{0,0} & \cdots & I_{0,n} \endvmatrix.$$
Thus, since $0 \leq i \leq n-1$, the top row will be equal to the $(n-i+1)th$ row, again making the quasideterminant 0 (c.f. {\cite{GGRW}}, prop. 1.4.6).

\vskip 12 pt
\textbf{Proof of Proposition:}\\
Let $n,m \in \mathbb{N}$ so that $n \neq m$. Suppose $n<m$. Now $p_{n}(x) = \sum_{k=0}^{n} a_{k}x^{k}$ for some $a_{0},\cdots,a_{n} \in R$. Thus $<p_{n}, q_{m}> = \sum_{k=0}^{n} a_{k}<x^{k}, q_{m}>$. For all $0 \leq k \leq n$, $k<m$ so by the lemma, $<x^{k}, q_{m}> = 0$. Thus $<p_{n}, q_{m}> = 0$. The case for $n>m$ is similar.

\vskip 12 pt
\underline{Remark}:
We note here that we can recover the construction of orthogonal polynomials in {\cite{GKLLRT}} from the construction above. Let R be the free associative algebra on generators $S_{0},S_{1},\cdots$ with  $S_{a+b} = I_{a,b}$ for all $a,b \in \mathbb{N}$. Following the notation of {\cite{GKLLRT}} let * be the anti-automorphism so $(S_{k})^{*} = S_{k}$ and $(\sum c_{i}x^{i})^{*} = \sum (c_{i})^{*}x^{i}$. A little examination shows that $q_{n}=p_{n}^{*}$. Thus $<p_{n},q_{m}> = <p_{n}, p_{m}>_{*}$, i.e. the collection $\{p_{n}\}$ is orthogonal with respect to the (very similar) inner product $<\cdot, \cdot>_{*}$ where $<\sum c_{i}x^{i}, \sum d_{j}y^{j}>_{*} = \sum c_{i}S_{i+j}(d_{j})^{*}$.

\vskip 12 pt

\section{Banded Matrices:}

For $i,j \in \mathbb{Z}_{\geq 0}$, let $E_{i,j}$ denote the matrix with rows and columns indexed by $\mathbb{Z}_{\geq 0}$ so that the $(i,j)$ entry is 1 and all other entries are 0. Let $a \leq 0$ and $b \geq 0$. $M_{[a,b]}$ is defined to be $span\{E_{i,j} : a \leq i-j \leq b\}$. We will refer to these matrices as ``banded''.
For example, the set of diagonal matrices is $M_{[0,0]}$.
Let $X \in M_{[a,b]}$ and $Y \in M_{[c,d]}$. 

\vskip 12 pt
\textbf{Lemma:} $X + Y \in M_{[min(a,c),max(b,d)]}$ and $XY \in M_{[a+c,b+d]}$.

\vskip 12 pt
\textbf{Proof:}\\
The proof that $X + Y \in M_{[min(a,c),max(b,d)]}$ is trivial.
Suppose $[XY]_{u,v} \neq 0$. Then $[X]_{u,w} \neq 0$ and $[Y]_{w,v} \neq 0$ for some w.
This implies $a \leq w-u \leq b$ and $c \leq v-w \leq d$. Adding these equations shows that $a+c \leq v-u \leq b+d$. Thus $XY \in M_{[a+c,b+d]}.$

\vskip 12 pt

\section{Recurrence Relations:}

In the commutative case, Bertola, Gekhtman and Szmigielski {\cite{BGS}} obtain a 4 term recurrence relation when $I_{a+1,b} + I_{a,b+1} = \alpha_{a}\beta_{b}$. This means there is a formula for $p_{n+1}$ in terms of $p_{n}, p_{n-1}$, and $p_{n-2}$ and a similar formula for $q_{n+1}$.$K$ is called the kernel of a system of biorthogonal polynomials if $<a(x),b(y)>$ = $\int\int a(x)b(y)K(x,y)dxdy$. The condition above corresponds to what the authors of this paper called the ``Cauchy kernel": $K(x,y) = \dfrac{1}{x+y}$. Below, we achieve similar, but longer, recurrences that correspond to kernels of the form $\dfrac{1}{f(x)+g(y)}$ where f and g are polynomials.

For all $n \in \mathbb{N}$, let $$p_{n} = \vmatrix I_{n,0} & \cdots & $\fbox{$I_{n,n}$} $\\ \vdots & \ddots & \vdots \\ I_{0,0} & \cdots & I_{0,n} \endvmatrix^{-1} \cdot \vmatrix I_{n,0} & \cdots & I_{n,n-1} & $\fbox{$x^{n}$} $\\ \vdots & \ddots & \vdots & \vdots \\ I_{0,0} & \cdots & I_{0,n-1} & 1 \endvmatrix$$ and let $$q_{n} = \vmatrix 1 & \cdots & $\fbox{$y^{n}$} $\\ I_{n-1,0} & \cdots & I_{n-1,n} \\ \vdots & \ddots & \vdots \\ I_{0,0} & \cdots & I_{0,n} \endvmatrix.$$
These are scalar multiples of the polynomials constructed in Section 2. Therefore they are biorthogonal. A quick check will show that we also have that $<p_{n},q_{n}> = 1$ for all $n \in \mathbb{N}$. Thus this system of polynomials is \underline{biorthonormal}.

\vskip 12 pt
\textbf{Theorem:}
Let $\{p_{k}\},\{q_{k}\}$ be any biorthonormal polynomials with bimoments $I$. Suppose there exist polynomials over the center of R $f(x) = \sum_{i=0}^n a_{i}x^{i}$ and $g(y) = \sum_{j=0}^m y^{j}b_{j}$ so that $\sum_{i=0}^n a_{i}I_{r+i,s} + \sum_{j=0}^m I_{r,s+j}b_{j} = \alpha_{r}\beta_{s}$ for all $r,s \in \mathbb{N}$. Then there exist $n+m+2$ term recurrence relations for $p_{i}$ and $q_{i}$. That is, we can express $p_{i+1}$ in terms of $p_{i}, \cdots, p_{i-n-m-2}$ and $q_{i+1}$ in terms of $q_{i}, \cdots, q_{i-n,m-2}$. The recurrences we achieve for $p_{i+1}$ and $q_{i+1}$ have polynomial coefficients for $p_{i}, p_{i-1}, q_{i}$, and $q_{i-1}$ and scalar coefficients for all other terms.

\vskip 12 pt
\textbf{Proof:}\\
Let $$\Lambda = \bmatrix 0 & 1 & 0 & \cdots \\ 0 & 0 & 1 & \ddots \\ \vdots & \vdots & \ddots & \ddots \endbmatrix.$$\\
Let $p(x)$ and $q(y)$ be column vectors with entries $p_{k}$ and $q_{k}$ respectively. Note that for each $k \in \mathbb{Z}_{\geq 0}$, $p_{k}$ and $q_{k}$ are polynomials of degree $k$ so for each so the products $p_{k}f(x)$ and $g(y)q_{k}$ can be written as a linear combination of $p_{n+k},\cdots,p_{1},p_{0}$ and $q_{m+k},\cdots,q_{1},q_{0}$ respectively.\\
Let $X$ and $Y$ be the infinite scalar matrices so that  $p(x)f(x) = Xp(x)$  and $g(y)q^{T}(y) = q^{T}(y)Y^{T}$. Since $<p(x),q^{T}(y)> = Id$, we know that $<p(x)f(x),q^{T}(y)> = X$ and $<p(x),g(y)q^{T}(y)> = Y^{T}$.\\
\vskip 9pt
Suppose $p_{k}(x) = \sum_{i=0}^k c_{i}x^{i}$ and $q_{l}(y) = \sum_{i=0}^l y^{i}d_{i}$. Let $\pi_{k} = \sum_{i=0}^k c_{i}\alpha_{i}$ and $\eta_{l} = \sum_{i=0}^l \beta_{i}d_{i}$.\\
$$(X+Y^{T})_{k,l} = <p_{k}(x)f(x),q_{l}(y)> + <p_{k}(x),g(y)q_{l}(y)> = $$
$$\sum_{i,j} c_{i}<f(x)x^{i},y^{j}>d_{j} + \sum_{i,j} c_{i}<x^{i},y^{j}g(y)>d_{j} = \sum_{i,j} c_{i}\alpha_{i}\beta_{j}d_{j} = \pi_{k}\eta_{l}.$$\\
If $\pi$ and $\eta$ are vectors with entries $\pi_{n}$ and $\eta_{n}$ respectively, then $X+Y^{T} = \pi\eta^{T} = D_{\pi}(\underline{1})(\underline{1}^{T})D_{\eta}$ where $D_{\pi}$ and $D_{\eta}$ are diagonal matrices with (i,i) entries $\pi_{i}$ and $\eta_{i}$, respectively.\\
\vskip 9pt
Let $A = (\Lambda - Id)D_{\pi}^{-1}X$ and $B^{T} = Y^{T}D_{\eta}^{-1}(\Lambda^{T} - Id)$. Since \underline{1} is a null vector of $\Lambda-Id$, $(\Lambda-Id)D_{\pi}^{-1}(X+Y^{T})=0$ and $(X+Y^{T})D_{\eta}^{-1}(\Lambda^{T}-Id)=0$. Then $A=-(\Lambda - Id)D_{\pi}^{-1}Y^{T}$ and $B^{T}=-X D_{\eta}^{-1}(\Lambda^{T} - Id)$.\\
\vskip 9 pt
We claim that A and B are banded matrices. Note that $X \in M_{[-\infty,n]}$ since $X_{i,j} = <p_{i}(x)f(x),q_{j}(y)> = 0$ if $i+n<j$ (because the degree $p_{i}*f(x)$ is less than the degree of $q_{j}$) and that $Y^{T} \in M_{[-m,\infty]}$ since $Y^{T}_{i,j} = <p_{i}(x),g(y)q_{j}(y)> = 0$ if $i>m+j$. Note also that $(\Lambda - I) \in M_{[0,1]}$.\\
\vskip 9pt
Applying the results we obtained for banded matrices, we see that $A = (\Lambda - Id)D_{\pi}^{-1}X \in M_{[-\infty,n+1]}$ and that $A = -(\Lambda - Id)D_{\pi}^{-1}Y^{T} \in M_{[-m,\infty]}$. Thus $A \in M_{[-m,n+1]}$. Similarly, $B^{T} \in M_{[-\infty,m+1]}$ and $B^{T} \in M_{[-n,\infty]}$ so $B^{T} \in M_{[-n,m+1]}$ and $B \in M_{[-m-1,n]}$.\\
\vskip 20pt
Recall that $p(x)f(x) = Xp(x)$  and $g(y)q^{T}(y) = q^{T}(y)Y^{T}$. Then $(\Lambda-Id)D_{\pi}^{-1}p(x)f(x)=(\Lambda-Id)D_{\pi}^{-1}Xp(x)=Ap(x)$ and $g(y)q^{T}(y)D_{\eta}^{-1}(\Lambda^{T} - Id)=
q^{T}(y)Y^{T}D_{\eta}^{-1}(\Lambda^{T} - Id)=q^{T}(y)B^{T}$.\\
\vskip 9pt
Thus examining the k-1th row of these equations gives the following n+m+2 term recurrence relations, as desired:
$$(\pi_{k}^{-1}p_{k} - \pi_{k-1}^{-1}p_{k-1})f(x) = \sum_{i=k-m}^{k+n+1} A_{k-1,i}p_{i},$$
$$g(y)(\eta_{k}^{-1}q_{k} - \eta_{k-1}^{-1}q_{k-1}) = \sum_{i=k-n}^{k+m+1} B_{k-1,i}q_{i}.$$\\
\vskip 12 pt

\section{Biorthogonal Analogue of Favard's Theorem:}

Favard’s theorem states that if $\{p_n(x)\}$ is a sequence of polynomials which obeys the usual 3-term recurrence relation then there exists an inner product for which these polynomials are orthogonal. Here we show that any two sequences of polynomials are biorthogonal with respect to some function, for which we construct the bimoments. It is important to note that no recurrence relation is required here.

\textbf{Theorem:}
Let $\{p_{n}\}, \{q_{n}\}$ be any set of polynomials over any division ring $R$ so that $p_{n}$ and $q_{n}$ are of degree $n$ for all $n \in \mathbb{N}$. For any $\{c_k\}_{k \in \mathbb{Z}_{\geq 0}}$ in $R$, there exists a unique set of bimoments for which $\{p_{n}\}, \{q_{n}\}$ is a biorthogonal system of polynomials and $<p_k,q_k>=c_k$.

\vskip 12 pt
\textbf{Proof:}\\
It is equivalent to show that there is a set of bimoments so that for all $a,b \in \mathbb{N}$, the following conditions hold:\\*
1) If $a<b$ then $<x^{a},q_{b}(y)>=0.$\\*
2) If $a>b$ then $<p_{a}(x),y^{b}>=0.$\\*
3) If $a=b$, then $<p_{a}(x),q_{b}(y)>=c_a.$\\*
We will define $I_{a,b}$ inductively on $a+b$. It is pivotal to note that the equations  $<x^{a},q_{b}(y)>=0$, $<p_{a}(x),y^{b}>=0$, and $<p_{a}(x),q_{b}(y)>=c_a$ do not involve bimoments of the form $I_{i,j}$ where $i+j>a+b$.
Recall that $p_0,q_0 \in R$. Let $I_{0,0}=p_{0}^{-1} c_0 q_{0}^{-1}$. Then $<p_0,q_0> = p_0 I_{0,0} q_0$ = 1 as desired.\\
Let $n \geq 1$ and  suppose for all a,b such that $a+b<n$, we have defined $I_{a,b}$ to satisfy the previous conditions.
For each $0 \leq i \leq n$ define $I_{i,n-i}$ as follows:
\vskip 6pt
\underline{Case 1:} If $i<n-i$ then the equation $<x^{i},q_{n-i}>$=0 is a linear equation whose variables (the bimoments) have all been defined except for $I_{i,n-i}$ due to the order in which the $I_{a,b}$'s are defined. Therefore there is a unique solution which we must define $I_{i,n-i}$ to be.
\vskip 6pt
\underline{Case 2:} Similarly, if $i>n-i$, the equation $<p_{i},y^{n-i}>$=0 has only one unknown and thus has a unique solution which we define $I_{i,n-i}$ to be.
\vskip 6pt
\underline{Case 3:} If $i=n-i$ then, again, the equation $<p_{i},q_{n-i}>=c_i$ has one unknown and we define $I_{i,n-i}$ to be the unique solution to this linear equation.\\
At each step we satisfy all the necessary conditions and have no choice so the bimoments constructed are the unique set for which $\{p_{n}\}, \{q_{n}\}$ is a biorthogonal system with $<p_k,q_k>=c_k$.

\end{document}